\documentclass[11pt]{amsart}

\usepackage[english]{babel}
\usepackage{csquotes}

\usepackage{amssymb}
\usepackage{amsthm}
\usepackage{amsmath}
\usepackage{mathrsfs}

\usepackage[a4paper, margin=3cm]{geometry}
\usepackage{color}
\usepackage{hyperref}
\usepackage{numprint}

\usepackage{enumerate}

\theoremstyle{plain}

\theoremstyle{definition}

\theoremstyle{remark}

\numberwithin{equation}{section}

\usepackage{tikz}
\usetikzlibrary{matrix,arrows}

\usepackage{mathrsfs}
\usepackage{amsmath, amsthm, amsfonts}
\usepackage[percent]{overpic}
\usetikzlibrary{matrix,arrows}

\usepackage{xypic}
\usepackage{color}
\usepackage{MnSymbol}
\usepackage{fullpage}
\usepackage{bbm}

\theoremstyle{definition}

\theoremstyle{remark}

\newcommand{\QQ}{\mathbb{Q}}
\newcommand{\EE}{\mathbb{E}}

\newcommand{\HH}{\mathcal{H}}
\newcommand{\PPP}{\mathcal{P}}

\newcommand{\RR}{\mathbb{R}}

\newcommand{\OO}{\mathcal{O}}

\newcommand{\ZZ}{\mathbb{Z}}

\newcommand{\PP}{\mathbb{P}}

\newcommand{\GL}{\text{GL}}
\newcommand{\dJ}{\text{dJ}}

\newcommand{\Pic}{\text{Pic}}
\newcommand{\Mg}{\mathcal{M}_g}
\newcommand{\Mgn}{\mathcal{M}_{g,n}}
\newcommand{\Mgb}{\overline{\mathcal{M}}_g}

\newcommand{\Cg}{\mathcal{C}_g}

\newcommand{\Jac}{\text{Jac}}
\newcommand{\rk}{\text{rk}}

\newcommand{\SL}{\mbox{SL}}

\newtheorem{Thm*}{Theorem*}
\theoremstyle{definition}

\title{Divisorial strata of abelian differentials}
\author{Scott Mullane\\
}
\date{\today}

\begin{document}
\thispagestyle{empty}

\maketitle

\begin{abstract}
We compute a closed formula for the class of the closure of the locus of curves in $\Mgb$ that admit an abelian differential of signature $\kappa=(k_1,...,k_{g-2})$.
\end{abstract}

\section{Introduction}
The strata of abelian differentials $\HH(\kappa)$ where $\kappa=(k_1,...,k_{n})$ with $k_i\in\ZZ_{>0}$ and $\sum k_i=2g-2$ is defined as all pairs $(C,\omega)$ known as translation surfaces where $C$ is a smooth genus $g$ curve and $\omega$ is a holomorphic differential on $C$ such that the zeros of $\omega$ have orders of type $\kappa$. Translation surfaces can be defined in an elementary way as polygons in the plane with certain side identifications and for this reason dynamics and rational billiards have provided much motivation for study in the area. 

Kontsevich and Zorich \cite{KontsevichZorich} have completely classified the connected components of $\HH(\kappa)$ for any $\kappa$. The natural $\GL^+(2,\RR)$ action on the plane induces an action on strata $\HH(\kappa)$. Recent breakthroughs by Eskin, Mirzakhani and Mohammadi \cite{EskinMirzakhani}\cite{EskinMirzakhaniMohammadi} have shown that the $\GL^+(2,\RR)$ orbit closure of any translation surface is a manifold in $\HH(\kappa)$ locally defined by linear equations of period coordinates with coefficients in $\RR$ and zero constant term. 

For any fixed $\kappa$, the closure of curves $C$ such that $(C,\omega)\in \HH(\kappa)$ form a subvariety in $\Mg$ and the Deligne-Mumford compactification $\Mgb$. If the length of $\kappa$ is $n=g-2$ then the resulting subvariety is of codimension one. Studying the algebraic subvarieties of a moduli space is an important aspect in the understanding of the birational geometry of a moduli space. Codimension one subvarieties are known as divisors and as $\Mg$ and $\Mgb$ are smooth, the group generated over $\QQ$ by linear equivalence classes of divisors take on a special importance as this space is isomorphic to $\Pic(\Mgb)\otimes \QQ$, the isomorphism classes of line bundles on $\Mgb$ modulo torsion. Harris and Mumford \cite{HarrisMumford}, Eisenbud and Harris \cite{EisenbudHarrisKodaira}, Farkas \cite{Farkas23}\cite{FarkasKoszul}\cite{Farkas}, Farkas and Verra \cite{FarkasVerra} have all used geometrically defined divisors to investigate the Kodaira dimension, structure of the Picard group and other aspects of various moduli spaces of curves.  

We calculate the class in $\Pic(\Mgb)\otimes\QQ$ of all such divisors coming from the strata $\HH(\kappa)$ where $\kappa$ has length $g-2$ which we denote $D_\kappa$. Where Kontsevich and Zorich \cite{KontsevichZorich} have shown there is more than one connected component we calculate the class of each connected component. This generalises the work of several mathematicians. Diaz \cite{Diaz} calculated the class of the closure of the locus of curves with an exceptional Weierstrass point in $\Mgb$. In our context this represents the class of $D_\kappa$ in $\Mgb$ for $\kappa=(g+1,1^{g-3})$. When $g=3$ this locus has two components and Cukierman \cite{Cukierman} \S5 calculated the class of the closure of the locus of hyperflexes on plane quartics which in our context is the odd spin structure component $D_\kappa^\text{odd}$ in $\overline{\mathcal{M}}_3$ for $\kappa=(4)$. Teixidor i Bigas \cite{Teixidor} calculated the divisor of curves with vanishing theta-null corresponding in our context to $D_\kappa^\text{even}$ in $\Mgb$ for $\kappa=(4,2^{g-3})$. Farkas and Verra \cite{Farkas}\cite{FarkasVerra} showed when the moduli space of even and odd spin curves has maximal Kodaira dimension by constructing effective divisors that when pushed down to $\Mgb$ correspond to the even and odd spin structure components $D_\kappa^\text{even}$ and $D_\kappa^\text{odd}$ in $\Mgb$ for $\kappa=(4,2^{g-3})$.

For $\kappa=(k_1,...,k_{g-2})$ we calculate  the class of the divisor $D_\kappa$ in $\Pic(\Mgb)\otimes \QQ$ to be
\begin{equation*}
D_\kappa=c_\lambda \lambda+\sum_{i=0}^{[g/2]}c_i\delta_i,
\end{equation*}
where
\begin{eqnarray*}
c_\lambda&=& \frac{7-g}{2(g-1)(g-2)}\left(  4(g-1)!\prod k_i^2+\sum_{k_i\geq 3}(k_i^2-k_i)\dJ[g-1;k_1,...,k_i-2,...,k_{g-2}]  \right)\\
&&+\frac{6}{g-1}\left((g-1)!\prod k_i^2+\sum_{i=1}^{g-2}k_i\left(\sum_{j=1}^{k_i-1} \dJ[g-1;j-1,k_1,...,k_{i-1},k_{i+1},...,k_{g-2}]     \right)\right),\\
c_0&=&\frac{-1}{2(g-1)}\biggl( (g-1)!\prod k_i^2+\sum_{i=1}^{g-2}k_i\left(\sum_{j=1}^{k_i-1} \dJ[g-1;j-1,k_1,...,k_{i-1},k_{i+1},...,k_{g-2}]     \right)\\
&&+ \frac{1}{2(g-2)}\biggl(  4(g-1)!\prod k_i^2+\sum_{k_i\geq 3}(k_i^2-k_i)\dJ[g-1;k_1,...,k_i-2,...,k_{g-2}]  \biggr)                                      \biggr),\\
c_{i}&=& \frac{-1}{2(g-i)-2}\biggl(\sum_{|I|=i,||I||\geq2i+1}\dJ[g-i;||I||-2i,k_j\text{ for $j\in I^C$}]\\
&&\left(i!\prod_{i\in I}k_i^2-\sum_{j\in I}(k_j+||I^C||-2(g-i)+1)\dJ[i;k_i\text{ for $i\in I-\{j\}$}]  \right)\\
&&+\sum_{|I|=i-1, ||I||\leq 2i-2} \dJ[i;k_j\text{ for $j\in I$}]\biggl((g-i)!(||I||-2i)^2\prod_{j\in I^C}k_j^2\\
&& -\sum_{\tiny{\begin{array}{cc}j\in I^C,\\k_j\geq2i-||I||+1\end{array}}}(k_j+(||I||-2i)+1) \dJ[g-i;k_j+(||I||-2i),k_i \text{ for $i\in I^C-\{j\}$}]    \biggr)  \biggr),  
\end{eqnarray*}
for $i=1,...,[g/2]$, where all formulas used are defined in \S\ref{sec:dJ}. In \S\ref{sec:comp} we calculate the class of each of the irreducible components when there is more than one irreducible component of the strata  providing results that agree with Teixidor i Bigas' \cite{Teixidor} result and Farkas and Verra's \cite{Farkas}\cite{FarkasVerra} divisor classes obtained on Cornalba's compactified spin moduli space. Our calculations rather take place in the moduli space of curves making use of what we know about the degeneration of abelian differentials and the theory of admissible covers. Computations from low genus seem to suggest that asymptotically, all such divisors have slope between $8$ and $9$.

To apply the method of test curves to our situation we must first introduce some relevant tools. We need to understand how the system of canonical divisors on a smooth curve degenerates as that curve degenerates to a singular stable curve of different types. Where applicable, we also need to understand how the notion of the spin structure of a canonical divisor degenerates as the underlying curve degenerates. Finally, we require methods for enumerating such occurrences on certain nodal curves which requires us to enumerate special holomorphic and meromorphic sections within linear equivalence classes of divisors on smooth curves. With all these tools in hand, we calculate the class of all divisorial strata of abelian differentials. 
\\
\\
\textbf{Acknowledgements.} I am very grateful to my advisor Dawei Chen for his guidance and many helpful discussions and comments during this project that will form part of my PhD thesis.

\section{Preliminaries}

\subsection{Divisor theory on $\Mgb$}
Let $\pi:\Cg\longrightarrow \Mgb$ be the universal curve. The \emph{Hodge bundle} $\EE$ on $\Mgb$ is defined as $\pi_*\omega$ where $\omega$ is the relative dualising sheaf of $\pi$. Hence $\EE$ is a vector bundle of rank $g$ and we define the \emph{Hodge class} as
\begin{equation*}
\lambda=c_1(\EE)=\det(\EE).
\end{equation*}
It is worth remarking that in general $\Mgb$ is not a fine moduli space and hence this universal curve may only exist up to finite base change. This will be of less significance to us as we consider the Picard group over $\QQ$. 

The Hodge class is an extension of the class defined on $\Mg$ and $\lambda$ generates $\Pic(\Mg)\otimes\QQ$, however, $\Pic(\Mgb)\otimes\QQ$ contains more classes.

The boundary $\Delta=\Mgb-\Mg$ of $\Mgb$ parameterising stable curves of genus $g$ with at least one node is codimension one. It is made up of components $\Delta_0,...,\Delta_{[g/2]}$, where $\Delta_0$ is the closure of the locus of stable curves that have a non-separating node and $\Delta_i$ is the closure of the locus of stable curves with a separating node that separates the curve into components with arithmetic genus $i$ and $g-i$ for $1\leq i\leq [g/2]$. These boundary divisors can intersect each other and self-intersect. Each $\Delta_i$ is irreducible and their classes in $\Pic(\Mgb)\otimes\QQ$ are denoted $\delta_i$. See \cite{AC,HarrisMorrison} for a more information.

For $g\geq 3$ the divisor classes $\lambda,\delta_0,...,\delta_{[g/2]}$ freely generate $\Pic(\Mgb)\otimes \QQ$. For $g=2$ the classes $\lambda,\delta_0$ and $\delta_1$ generate $\Pic(\overline{\mathcal{M}}_2)\otimes \QQ$ with the relation
\begin{equation*}
\lambda=\frac{1}{10}\delta_0+\frac{1}{5}\delta_1.
\end{equation*}

\subsection{Strata of abelian differentials}
A \emph{signature} $\kappa=(k_1,...,k_n)$ is a partition of $2g-2$ with all $k_i\in \ZZ_{>0}$. We define the \emph{stratum of abelian differentials with signature $\kappa$} as
\begin{equation*}
\HH(\kappa):=\{ (C,\omega) \hspace{0.15cm}|\hspace{0.15cm} g(C)=g,\hspace{0.05cm} \omega\in H^0(C,K_C)\text{ such that }(\omega)_0=k_1p_1+...+k_np_n,\text{ for $p_i$ distinct}\}
\end{equation*}
that is, the space of abelian differentials with prescribed multiplicities of zeros given by $\kappa$. By relative period coordinates $\HH(\kappa)$ has dimension $2g+n-1$. One advantage of this construction is that an abelian differential $(C,\omega)$ can be represented by polygons in the plane with certain side identifications. The deformation space of such an abelian differential can then be visualised very concretely as perturbing the sides of the polygons. Almost paradoxically, this visual tool of the deformation space is not available when we forget the differential $\omega$ and only consider the deformations of the genus $g$ curve $C$.

Survey articles \cite{Wright} and \cite{Zorich} provide broad introductions to the theory of the strata of abelian differentials. 

A related object is the \emph{stratum of canonical divisors with signature $\kappa$} which we define as
\begin{equation*}
\PPP(\kappa):=\{[C,p_1,...,p_n]\in \Mgn   \hspace{0.15cm}| \hspace{0.15cm}k_1p_1+...+k_np_n\sim K_C\}.
\end{equation*}
Observe that the zeros of the canonical divisor with $k_i=k_j$ are ordered. Forgetting this ordering of the zeros we obtain a finite cover of the projectivisation of $\HH(\kappa)$. If all $k_i$ are distinct then  $\PPP(\kappa)$ is isomorphic to the projectivisation of $\HH(\kappa)$. Hence we have the dimension of $\PPP(\kappa)$ is $2g+n-2$.

A \emph{theta characteristic} on a smooth curve $C$ is a line bundle $\eta$ on $C$ such that $\eta^{\otimes 2}\sim K_C$. The parity of $h^0(C,\eta)$ is known as a \emph{spin structure} of $\eta$ and Mumford \cite{Mumford} showed that this parity is deformation invariant. Consider an abelian differential $(C,\omega)$ where $\omega$ has signature $\kappa=(2k_1,..,2k_n)$. Then there is a natural choice of theta characteristic for this abelian differential
\begin{equation*}
\eta\sim \sum_{i=1}^n k_ip_i
\end{equation*}
and the loci $\HH(\kappa)$ and $\PPP(\kappa)$ are reducible and break up into disjoint components with even and odd parity of $h^0(C,\eta)$. A \emph{hyperelliptic differential} of type $\kappa$ is a differential on a hyperelliptic curve with the minimum number of zeros occurring at ramification points of the hyperelliptic involution known as Weierstrass points. The subvariety of hyperelliptic differentials in $\HH(\kappa)$ has dimension $2g+(n-m)/2$ where $m$ is the number of zeros that occur at Weierstrass points in each hyperelliptic differential. Kontsevich and Zorich \cite{KontsevichZorich} showed that there can be at most $3$ connected components in total of $\HH(\kappa)$ and hence $\PPP(\kappa)$, corresponding to the case that the hyperelliptic differentials become a connected component of $\HH(\kappa)$ distinct from the remaining differentials that provide two further connected components based on odd or even spin structure.

\subsection{Limit linear series and degeneration of canonical divisors}
Einsenbud and Harris \cite{EisenbudHarrisLimit} develop the theory of limit linear series. A \emph{linear series} of degree $d$ and dimension $r+1$, or $g^r_d$ on $C$ is an $(r+1)$-dimensional vector space of linearly equivalent divisors on $C$. Hence we can express a $g^r_d$ on a curve $C$ as a pair $(L,V)$ where $L$ is a line bundle on $C$ and $V$ is an $(r+1)$-dimensional subspace of $H^0(C,L)$.

For any linear series $V$ on a curve $C$ and any smooth point $p\in C$ we have the 
\emph{vanishing sequence}
\begin{equation*}
0\leq a_0(V,p)<a_1(V,p)<\cdot\cdot\cdot<a_r(V,p)
\end{equation*}
which are just the orders of vanishing at $p$ of the elements of the $r$-dimensional $\PP V$. The \emph{ramification sequence} $b_i(V,p)$ is defined by
\begin{equation*}
b_i(V,p)=a_i(V,p)-i
\end{equation*}
for $i=0,1,...,r$. The \emph{ramification index} of $V$ at $p$ is given by
\begin{equation*}
\beta(V,p)=\sum_{i=0}^rb_i(V,p).
\end{equation*}
The Pl\"ucker formula gives the ramification over all points. If $V$ is any $g_d^r$ on a smooth curve $C$ of genus $g$,
\begin{equation*}
\sum_{p\in C}\beta(V,p)=(r+1)d+(r+1)r(g-1).
\end{equation*}
Harris and Morrison \cite{HarrisMorrison} provide an introduction to linear series and the technique of limit linear series.

Every smooth curve $C$ of genus $g$ has a unique $g^{g-1}_{2g-2}$ known as the \emph{Canonical series} equal to $\PP H^0(C,K_C)$. The Pl\"ucker formula gives us the total ramification of this series as $(g+1)g(g-1)$. We call the points on a curve with non-trivial ramification \emph{Weierstrass points} and on a general curve there are $(g+1)g(g-1)$ distinct Weierstrass points each with simple ramification. We are interested in how the Canonical series degenerates as $C$ degenerates to a nodal curve.

A nodal curve $C$ is of \emph{compact type} if removing any node makes the curve disconnected. Eisenbud and Harris \cite{EisenbudHarrisLimit} introduced the theory of \emph{limit linear series} as the limiting objects of a family of $g^r_d$'s on smooth curves degenerating to a nodal curve of compact type. If $C$ is a curve of compact type with components $C_1,...,C_n$ then a (refined) limit $g^r_d$ is a collection $\{V_i\}$ of a $g^r_d$ on each $C_i$ such that if $C_i$ and $C_j$ intersect at a node $q$ then
\begin{equation*}
a_m(V_i,q)+a_{r-m}(V_j,q)=d
\end{equation*}   
for $m=0,...,r$. Eisenbud and Harris remarked that this method also applies to a larger class of curves. We define a curve of \emph{pseudocompact type} to be a curve in which every node is either a separating node (disconnects the curve) or a self node (the self intersection of an irreducible component in the normalisation of the curve).

A canonical divisor with signature $\kappa=(k_1,...,k_n)$ on a smooth curve $C$ is a section of the canonical series on $C$ with zeros of type $\kappa$. A \emph{limit canonical divisor} on a curve of compact type with signature $\kappa$ is a section of a limit canonical series where the zeros away from the nodes are of type $\kappa$. Eisenbud and Harris \cite{EisenbudHarrisWeierstrass} showed that a limit canonical divisor with signature $\kappa$ on a curve of compact type is the limit of canonical divisors of signature $\kappa$ on smooth curves and investigated the location of the limits of Weierstrass points on such curves of compact type. For a curve of pseudocompact type the only alteration that we must make is to allow zeros at the non-separating nodes. In the situation where a zero of $\kappa$ is occurring at a non-separating node we may want to consider this limit as the curve blown up at the node and a $\PP^1$-bridge inserted with a marked point on this bridge representing the zero. We will be clear which limit we are considering. Esteves and Medeiros \cite{EstevesMedeiros} investigated limit canonical divisors on curves with two components that were not of compact type while Osserman \cite{OssermanNoncompact} investigates limit linear series in general for curves not of compact type.

\subsection{Degeneration of theta characteristics and spin structures}
Cornalba \cite{Cornalba} discusses how theta characteristics degenerate to nodal curves including curves of pseudocompact type. We first consider a non separating node. Let $C$ be a curve with only one node which is non-separating. Let $\tilde{C}$ be the normalisation of $C$ with $x$ and $y$ the preimages of the node under the normalisation map. Cornalba  \cite{Cornalba} showed that there are two types of theta characteristics on such a curve. For the first type we begin with a line bundle $\tilde{\eta}$ on $\tilde{C}$ such that $\tilde{\eta}^{\otimes 2}\sim K_{\tilde{C}}+x+y$. We then observe that as $\tilde{\eta}^{\otimes 2}\sim K_{\tilde{C}}+x+y$ and $K_{\tilde{C}}+x$ has a base point for any $x$ we have a section of $H^0(\tilde{C},\tilde{\eta})$ has a zero at $x$ if and only if it has a zero at $y$. Hence the sections of $H^0(\tilde{C},\tilde{\eta})$ that vanish at $x$ and $y$ form a codimension one locus of $H^0(\tilde{C},\tilde{\eta})$. Having a line bundle $\tilde{\eta}$ on $\tilde{C}$ there are two ways that this bundle can descend to a bundle $\eta$ on $C$. We can glue sections as $f(x)=f(y)$ or $f(x)=-f(y)$ as in both cases the square of these will agree at the node. These two possibilities will have $h^0(C,\eta)$ differing by $+1$ and hence represent odd and even spin structures.

The second possibility is that we blow up at the node and insert an exceptional $\PP^1$-bridge between $x$ and $y$. Here we have the theta characteristics are
\begin{equation*}
(\tilde{\eta},\OO(1))
\end{equation*} 
where $\tilde{\eta}^{\otimes 2}\sim K_{\tilde{C}}$ and we consider the global sections to be glued together from $H^0(\tilde{C},\tilde{\eta})$ and $H^0(\PP^1,\OO(1))$ at the nodes $x$ and $y$. But $h^0(\PP^1,\OO(1))=2$ and hence the values at $x$ and $y$ completely determine the section on $\PP^1$. The parity of such a theta characteristic is thus $h^0(\tilde{C},\tilde{\eta})$ mod $2$.

We are now ready to consider theta characteristics on a curve of pseudocompact type. We consider a curve $C$ of pseudocompact type with irreducible components $C_1,...,C_k$. We first blow up at every separating node and insert an exceptional component i.e. a $\PP^1$ between the two components. Then a theta characteristic on the curve $C$ is
\begin{equation*}
(\eta_1,...,\eta_k,\{\OO(1)\}_{i=1}^{k-1})
\end{equation*}
where $\eta_i$ is a theta characteristic on $C_i$ and $\OO(1)$ is a line bundle of degree one on the exceptional $\PP^1$ components. This gives the total degree $\sum_{i=1}^k(g_i-1)+(k-1)=g-1$ as expected. Observing that $h^0(\PP^1,\OO(1))=2$ we see that the parity of this spin structure is given by
\begin{equation*}
\sum_{i=1}^k h^0(C_i,\eta_i) \text{   mod } 2
\end{equation*}
where if any component $C_i$ has self nodes then the $\eta_i$ is of the types discussed earlier. Dealing with spin structures can be subtle. For example, Chen and Gendron (\cite{ChenDegen}\S7, \cite{Gendron} \S7) both show that though $\PPP(4)^\text{odd}$ and $\PPP(4)^\text{hyp}=\PPP(4)^\text{even}$ are disjoint in $\mathcal{M}_{3,1}$, their closures intersect in $\overline{\mathcal{M}}_{3,1}$ by providing two families of smooth marked curves of different parities limiting to the same marked nodal curve.

\subsection{de Jonquieres' Formula}\label{sec:dJ}
The total number of sections of a general $g_d^r$ on a genus $g$ curve with ordered zeros of order $k_i$ for $i=1,...,\rho$ with $\sum k_i=d$ and $\rho=d-r$ is
\begin{equation*}
\dJ[g;k_1,...,k_\rho]=\frac{g!}{(g-\rho-1)!}\prod_{i=1}^{\rho}k_i\left( \sum_{j=0}^{\rho-1}\left(\frac{(-1)^j}{g-\rho+j}\sum_{|I|=j}\left(\prod_{i\nin I}k_i    \right)  \right)+\frac{(-1)^\rho}{g}\right).
\end{equation*}
where $I$ is a subset of $\{1,...,\rho\}$ and $|I|$ denotes the number of elements in $I$. For later convenience of notation we define 
\begin{equation*}
||I||:=\sum_{i\in I}k_i.
\end{equation*}
In the case that all $k_i$ are distinct, this is an equivalent formula to that presented in \cite{ACGH} on page $359$. Our version is more convenient for computational purposes and is only a slight variation of that developed in \cite{Coolidge}  page $288$ where we are interested in marking the zeros of the section. We will use the convention that $\dJ[1;\emptyset]=1$.

\subsection{The Picard variety method}
A situation that will often present itself is that we will want to know the number of solutions to a particular equation occurring in the Picard group of a general curve $C$ of genus $g$. In general, we will want to know for a general curve $C$ of genus $g$, how many $(p_1,...,p_g)\in C^g$ there are such that 
\begin{equation*}
\sum_{i=1}^gk_ip_i\sim L
\end{equation*}
where $L$ is a specified line bundle of order $d=\sum_{i=1}^gk_i$ and $k_i\ne0$, $k_i\in \ZZ$ with only a finite number of sections of the type required. We will follow the treatment of \cite{ChenTarasca} \S 2. We consider the map
\begin{eqnarray*}
\begin{array}{cccccc}
f:C^g&\longrightarrow &\Pic^d(C)\\
(p_1,...,p_g)&\longmapsto&\sum_{i=1}^gk_ip_i.
\end{array}
\end{eqnarray*}
The fibre of this map above $L\in\Pic^d(C)$ will give us precisely the solutions of interest. We observe that the domain and range of $f$ are both of dimension $g$. Hence once we have identified that there are only a finite number of solutions for our specific $k_i$, our answer is simply the degree of the map $F$. Take a general point $e\in C$ and consider the isomorphism
\begin{eqnarray*}
\begin{array}{cccccc}
h:\Pic^d(C)&\longrightarrow &J(C)\\
L&\longmapsto&L\otimes\OO_C(-de).
\end{array}
\end{eqnarray*}
Now let $F=h\circ f$. Then we have $\deg F=\deg f$. We observe
\begin{equation*}
F(p_1,...,p_g)=\OO_C\biggl(\sum_{i=1}^gk_i(p_i-e)\biggr).
\end{equation*}
Let $\Theta$ be the fundamental class of the theta divisor in $J(C)$. By \cite{ACGH} \S1.5 we have
\begin{equation*}
\deg \Theta^g=g!
\end{equation*}
and the locus of $\OO_C(k(x-e))$ for varying $x\in C$ has class $k^2\Theta$ in $J(C)$. Hence
\begin{eqnarray*}
\deg F&=&\deg F_*F^*([\OO_C])\\
&=&\deg \left(\prod_{i=1}^g k_i^2\Theta\right)\\
&=&g!\left(\prod_{i=1}^g k_i^2\right)
\end{eqnarray*}
In practice we may want to discount this number by any specific solutions that we may be omitting for some reason. For example, we will be omitting any solutions where $p_i=p_j$ for $i\ne j$. In this case we will need to know not only the existence of any specific solutions that we are discounting by, but also the multiplicity of these solutions. We calculate the multiplicity by investigating the branch locus of $F$. First we look locally analytically at $F$ around each point. If $f_0d\omega,...,f_{g-1}d\omega$ is a basis for $H^0(C,K_C)$, then locally analytically the map becomes
\begin{eqnarray*}
(p_1,...,p_g)&\longmapsto&\biggl(\sum_{i=1}^gk_i\int_e^{p_i}f_0d\omega,...,\sum_{i=1}^gk_i\int_e^{p_i}f_{g-1}d\omega\biggr)
\end{eqnarray*}
modulo $H_1(C,K_C)$. The map on tangent spaces at any fixed point $(p_1,...,p_g)\in C^g$ is the Jacobian of $F$ at the point, which is
\begin{equation*}
\Jac(F)_{(p_1,...,p_g)}=\text{diag}(k_1,...,k_g)\begin{pmatrix}f_0(p_1)&...&f_0(p_g)\\
f_1(p_1)&...&f_1(p_g)\\
...&...&...&\\
f_{g-1}(p_1)&...&f_{g-1}(p_g)   \end{pmatrix}
\end{equation*}
Ramification in the map $F$ occurs when the map on tangent spaces is not injective which takes place at the points where $\rk(\Jac(F))<g$. The ramification index at a point $(p_1,...,p_g)\in C^g$ will be equal to the vanishing order of the determinant of $\Jac(F)_{(p_1,...,p_g)}$ at the point. 

We observe that there are two components to the branch locus of $F$. 
\begin{eqnarray*}
\Delta&=& \{(p_1,...,p_g)\in C^g\hspace{0.15cm}|\hspace{0.15cm}p_i=p_j \text{ for some $i\ne j$}\}      \\
\mathcal{K}&=&\{(p_1,...,p_g)\in C^g\hspace{0.15cm}|\hspace{0.15cm}h^0(C,K_C-p_1-...-p_g)>0\}
\end{eqnarray*}
where $\mathcal{K}$ is irreducible and $\Delta$ has $g(g-1)/2$ irreducible components defined by
\begin{equation*}
\Delta_{i,j}= \{(p_1,...,p_g)\in C^g\hspace{0.15cm}|\hspace{0.15cm}p_i=p_j \}
\end{equation*}
for $i,j=1,...,g$ and $i< j$. Hence finding the order of any point in the branch locus will simply be a matter of investigating how these loci meet at the particular point.

\section{Test curves}
Before we start to compute the class of the locus $D_\kappa$ we must be precise in our definition of this locus. Consider the locus of pointed curves $[C,p_1,...,p_{g-2}]\in \mathcal{M}_{g,g-2}$ such that
\begin{equation*}
\sum_{i=1}^{g-2}k_ip_i\sim K_C.
\end{equation*}
Pushing this locus down to a codimension one locus in $\Mg$ and taking the closure we obtain the divisor $D_\kappa$ in $\Mgb$. 

We now allow 
\begin{equation*}
D_\kappa=c_\lambda\lambda+\sum_{i=0}^{[g/2]}c_i\delta_i
\end{equation*}
for unknown coefficients $c_\lambda,c_i$. Intersecting either side of this equation with a test curve we will produce a relation between the coefficients. With enough relations we can recover all of the coefficients. Harris and Morrison \cite{HarrisMorrison} \S3F provide an introduction to the method of test curves.

\subsection{Test curve $A$}
Take a pencil of plane cubics. Attach one base point to a genus $g-1$ curve $C$ at a general point $y$ on the curve. 

\begin{figure}[htbp]
\begin{center}
\begin{overpic}[width=0.35\textwidth]{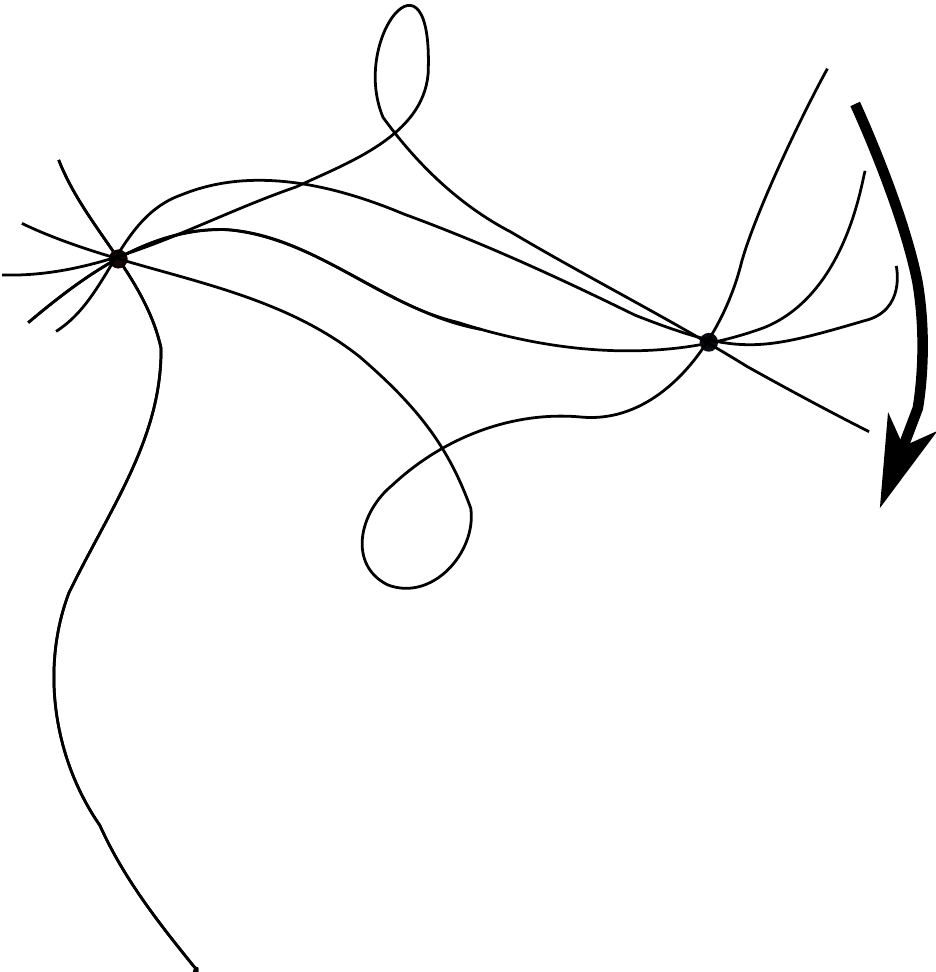}

\put (9,65){$y$}
\put (-35,18){$g(C)=g-1$}
\put (-5,35){$C$}
\put (60,35){A pencil of plane cubics}

\end{overpic}
\end{center}
\end{figure}

This is a standard test curve and it is well-known \cite{HarrisMorrison} \S3F that $A\cdot \lambda=1,A\cdot\delta_0=12,A\cdot\delta_1=-1$, giving
\begin{equation*}
A\cdot D_\kappa=c_\lambda+12c_0-c_1=0.
\end{equation*}
We know this intersection to be zero as if the attaching point on the genus $g-1$ curve is general then no limits of the type we are considering can occur for any genus $1$ curve attached to a genus $g-1$ curve at a general point. Any such limit canonical divisor would restrict on the genus $g-1$ component to be
\begin{equation*}
\sum_{i\in I}k_ip_i+(2g-4-\sum_{i\in I}k_i)y\sim K_C
\end{equation*}
for some subset $I\subseteq \{1,...,{g-2}\}$. However, for any choice of subset $I$ this would require either the curve $C$ or the point $y$ to be special, providing a contradiction.

\subsection{Test curve $B$}
Take a smooth general genus $g-1$ curve $C$. Create a node by identifying one non-special fixed point $y$ on the curve with another point $x$ that varies in the curve. 

\begin{figure}[htbp]
\begin{center}
\begin{overpic}[width=0.25\textwidth]{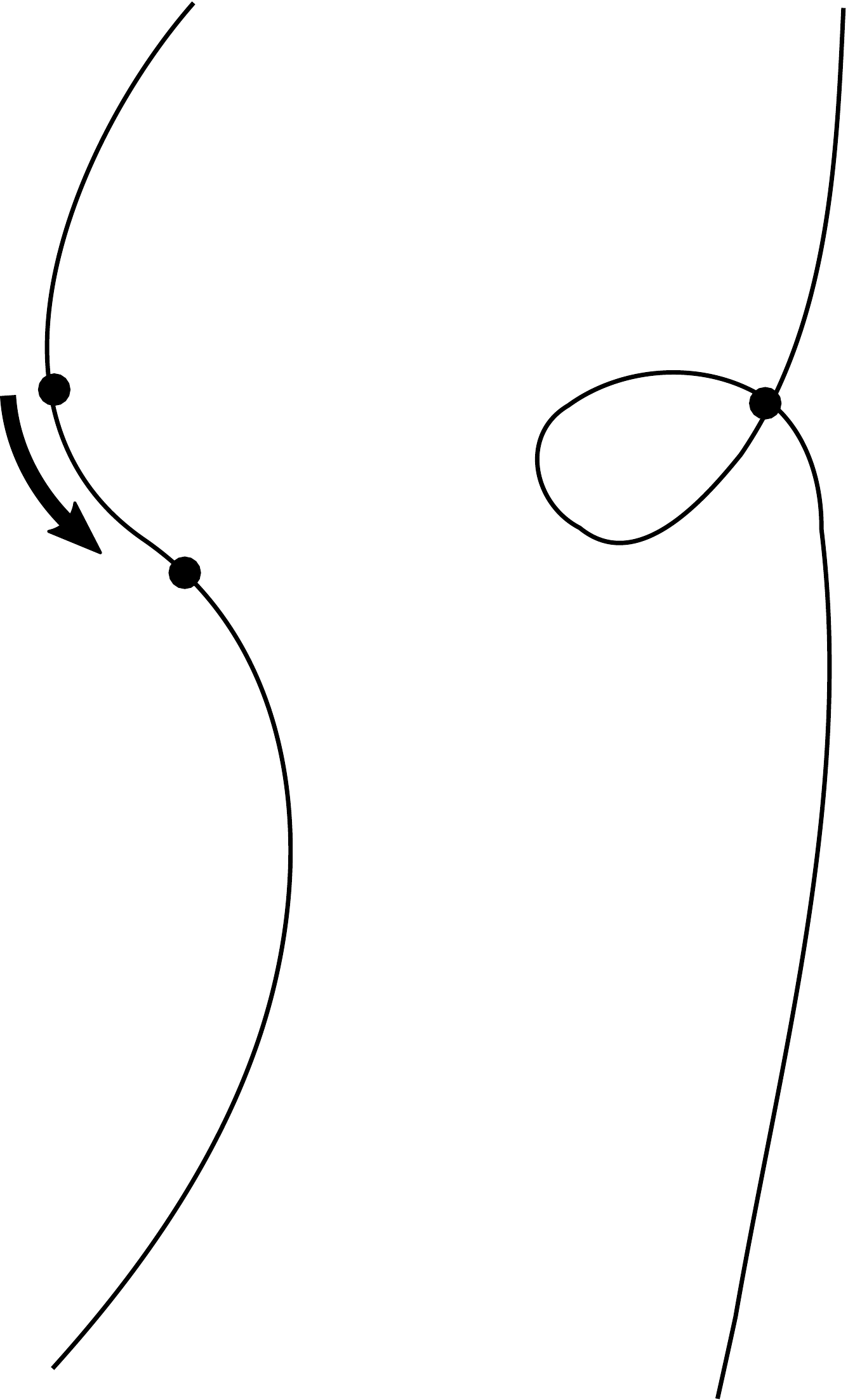}

\put (8, 71){$x$}
\put (18,60){$y$}
\put (-15,35){$g(C)=g-1$}
\put (7,18){$C$}
\put (60,68){$x\sim y$}

\end{overpic}
\end{center}
\end{figure}

This is again a well-known test curve \cite{HarrisMorrison} \S3F and $B\cdot \lambda=0$, $B\cdot \delta_0=2-2g$, $B\cdot \delta_{1}=1$ with the intersection with all other boundary components being zero, giving
\begin{equation*}
B\cdot D_\kappa=(2-2g)c_0+c_1.
\end{equation*}
To calculate this intersection directly we have that on a nodal curve of this type for $x\ne y$ the limits of differentials of type $\kappa$ will be of two types. Limits of the first type will satisfy
\begin{equation*}
\sum_{i=1}^{g-2}k_ip_i\sim K_C+x +y
\end{equation*}
for $p_i\ne x$ or $y$. Limits of the second type will satisfy
\begin{equation*}
jx+(k_m-j)y+\sum_{i\ne m}k_ip_i\sim K_C+x+y
\end{equation*}
for some $k_m$ with $j=1,...,k_m-1$. We will start by counting all limits of the first type. As $x$ is varying we consider the map
\begin{eqnarray*}
\begin{array}{cccc}
C^{g-1}&\rightarrow& \Pic^{2g-3}(C)\\
(p_1,...,p_{g-2},x)&\mapsto&\sum_{i=1}^{g-2}k_ip_i-x
\end{array}
\end{eqnarray*} 
which has degree $(\prod k_i^2)\Theta^2=(\prod k_i^2)(g-1)!$. But we must discount any solutions where $p_i=x$ or $y$. But as $K_C+x$ has a base point $x$ we see that there are no solutions with $p_i=x$ or $y$ as this would imply that $K_C$ has a special section which is not possible on a general curve $C$ or that $y$ forms part of a special section in a general curve $C$ which is also not possible as we chose $y$ to be general.

A limit of the second type will satisfy 
\begin{equation*}
(j-1)x+\sum_{i\ne m}k_ip_i\sim K_C-(k_m-j-1)y
\end{equation*}
for some $k_m$ with $j=1,...,k_m-1$. We observe immediately that there will be no solutions for $j=1$ as $C$ is a general curve and $y$ a general point. For other $j$ there are 
\begin{equation*}
\dJ[g-1;j-1,k_1,...,k_{m-1},k_{m+1},...,k_{g-2}]
\end{equation*}
solutions of this type and it only remains to find the order of such solutions. We would like to know in a general family of smooth curves which have a canonical divisor of type $\kappa$ specialising to this nodal curve, how many different canonical divisors are specialising to this canonical divisor with a zero of order $k_m$ at the node.

\begin{figure}[htbp]
\begin{center}
\begin{overpic}[width=0.9\textwidth]{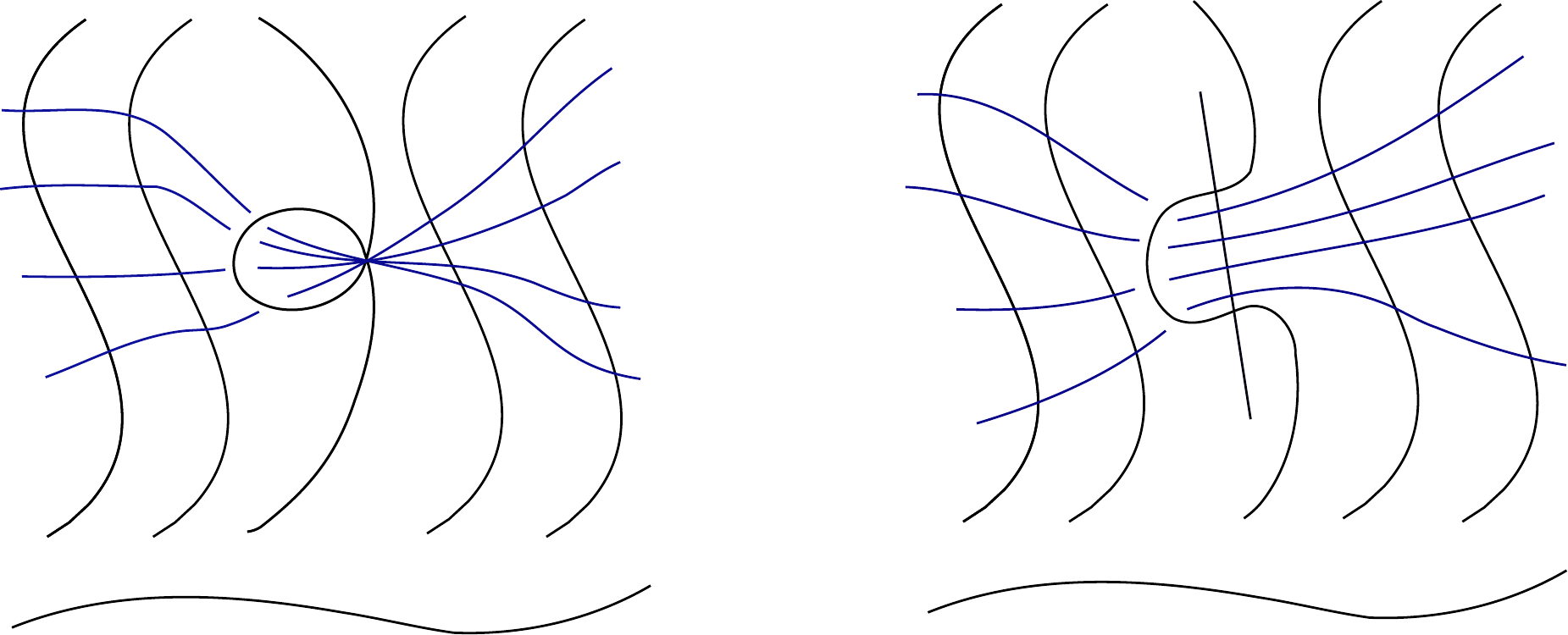}

\end{overpic}
\end{center}
\end{figure}

To find this order we observe that in this case we are actually considering a solution in the boundary of $\overline{\mathcal{M}}_{g,g-2}$ and hence we must blow up at the node to obtain this solution with a single zero of order $k_m$ sitting on a  $\PP^1$-bridge between $x$ and $y$. Osserman \cite{OssermanNoncompact} investigated limit linear series on curves of non-compact type and Esteves and Medeiros \cite{EstevesMedeiros} investigated limit canonical divisors on curves with two components. In our situation we have a so-called ``banana curve" with a smooth genus $g-1$ curve $C$ attached at points $x$ and $y$ to a $\PP^1$. We know the series $|K_C+x+y|$ is the $g^{g-1}_{2g-2}$ that appears in the limit $g^r_d$  on the component $C$. We have chosen $C$ to be a general curve and $y$ a general point. We also know that $x$ and $y$ sit in a special section in $|K_C+x+y|$. The difference between this case and the non-compact case is that the $g^{g-1}_{2g-2}$ that will appear in the limit on the $\PP^1$ component will adhere to gluing conditions at two nodes and how $x$ and $y$ are related in the series $|K_C+x+y|$ provides conditions on the $g^{g-1}_{2g-2}$ on our exceptional $\PP^1$. We observe that by imposing vanishing of $j$ at $x$ and $k_m-j$ at $y$ in a section in $|K_C+x+y|$ we impose total vanishing of $2g-2-k_m$ at the nodes in $\PP^1$ leaving us with a $g^{k_m-1}_{k_m}$ which has $k_m$ simple ramification points. Hence the order of our solution is $k_m$.

Putting this together we have
\begin{equation*}
B\cdot D_\kappa=(g-1)!\prod k_i^2+\sum_{i=1}^{g-2}k_i\left(\sum_{j=1}^{k_i-1} \dJ[g-1;j-1,k_1,...,k_{i-1},k_{i+1},...,k_{g-2}]     \right).
\end{equation*}

\subsection{Test curve $C_1$}
Let $Y$ be a genus $g-1$ curve and $X$ be an elliptic curve. Attach $X$ to $Y$ at a general point $x$ in $X$ and allow the attaching point $y$ in $Y$ to vary. 

\begin{figure}[htbp]
\begin{center}
\begin{overpic}[width=0.5\textwidth]{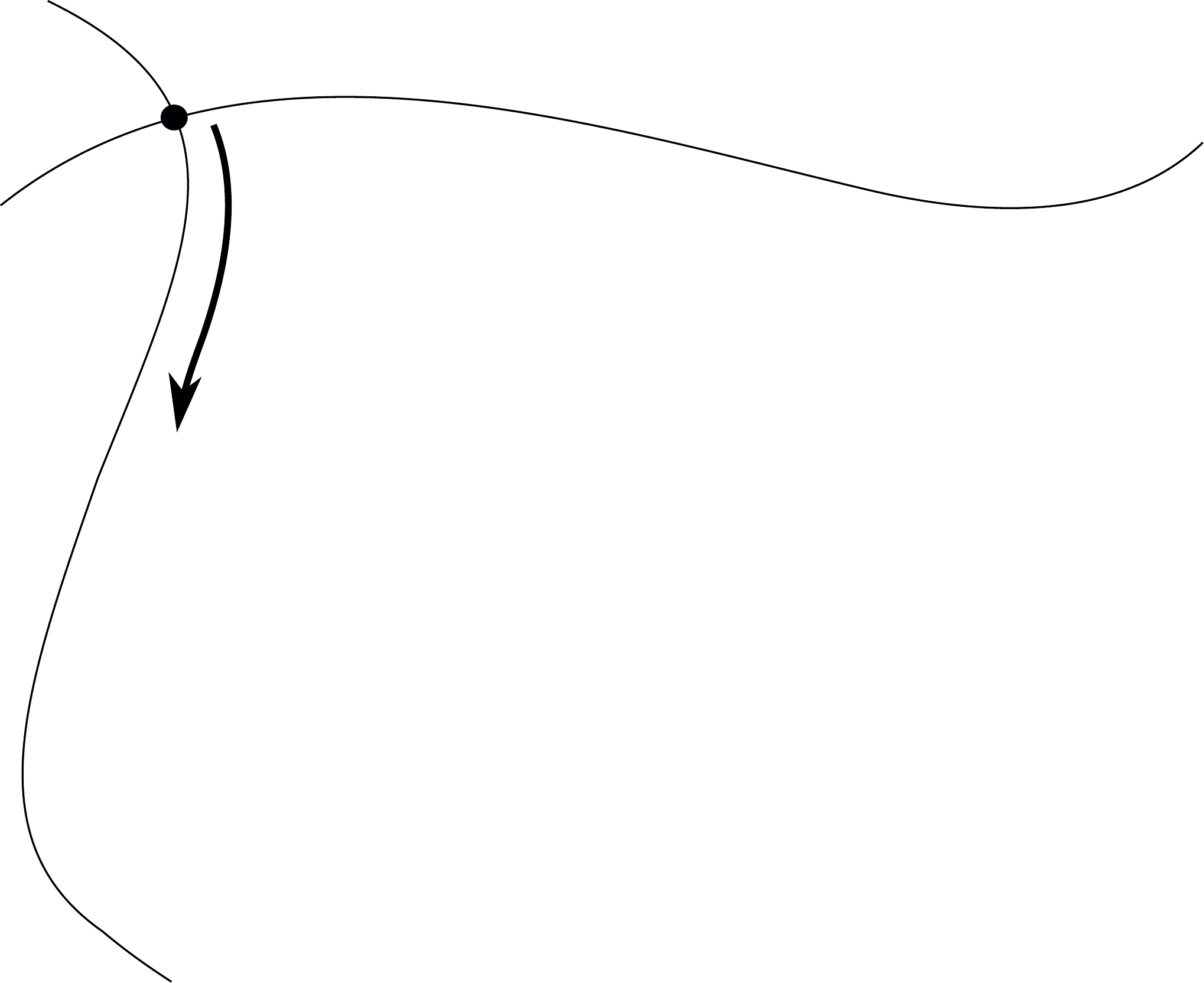}
\put (15, 74){$x$}
\put (65,72){$X$}

\put (45,63){$g(X)=1$}
\put (12,67){$y$}
\put (11,30){$g(Y)=g-1$}
\put (11,18){$Y$}
\end{overpic}
\end{center}
\end{figure}

We observe that $C_1\cdot D_\kappa=(2-2(g-1))c_{1}=(4-2g)c_1$. 
 
To calculate this intersection directly we must locate the limits of the type $\kappa$ in our test curve. If the limit of any $p_m$ occurs on $X$ we have that in the $Y$-aspect
\begin{equation*}
(k_m-2)y+\sum_{i\ne m}k_ip_i\sim K_Y.
\end{equation*}
But this is only possible on a general curve $Y$ if $k_m\geq 3$. In this case there are 
\begin{equation*}
\dJ[g-1;k_1,...,k_m-2,...,k_{g-2}]
\end{equation*} 
solutions. While on the $X$ component there are $k_m^2-1$ points satisfying $k_mp_m=k_mx$ with $p_m\ne x$.

The other possibility is that all $g-2$ points occur in the $Y$ component. in this case we have 
\begin{equation*}
\sum_{i=1}^{g-2}k_ip_i-2y\sim K_Y.
\end{equation*}
By the Picard variety method we know that there are $4(g-1)!\prod k_i^2$ solutions and we must discount for any solutions where $p_i=y$. But we just found solutions of this type. These solutions will each have order $k_i-1$ by the order of the intersection of the components of the ramification loci of the relevant Jacobi map at these points.

The last possibility is that there are no $p_i$ on $X$ and $g-3$ special points on $Y$ with the last special point  $p_m$ lying on a $\PP^1$-bridge between $x$ and $y$. Let $s$ be the point on $\PP^1$ glued to $x$ and $t$ the point glued to $y$. Then the limit canonical divisor of the type we require restricts on the $\PP^1$-bridge to be
\begin{equation*}
-k_mt+k_mp-2s\sim K_{\PP^1}.
\end{equation*}
As there is only one pole in the $Y$-aspect at $y$ the residue must be zero. Similarly the residue at $x$ is zero hence our section on $\PP^1$ must have zero residues at $s$ and $t$. We see that this is only possible for $k_m=0$. Refer to (\cite{EisenbudHarrisWeierstrass}, \S2) for a discussion of residues that includes this situation.
 
Hence we obtain
\begin{equation*}
C_1\cdot D_\kappa=4(g-1)!\prod k_i^2+\sum_{k_i\geq 3}(k_i^2-k_i)\dJ[g-1;k_1,...,k_i-2,...,k_{g-2}].
\end{equation*}

\subsection{Tests curve $C_i$}
Let $Y$ be a genus $g-i$ curve and $X$ be a genus $i$ curve for $1\leq i\leq g-2$. Attach $X$ to $Y$ at a general point $x$ in $X$ and allow the attaching point $y$ in $Y$ to vary. 
\begin{figure}[htbp]
\begin{center}
\begin{overpic}[width=0.5\textwidth]{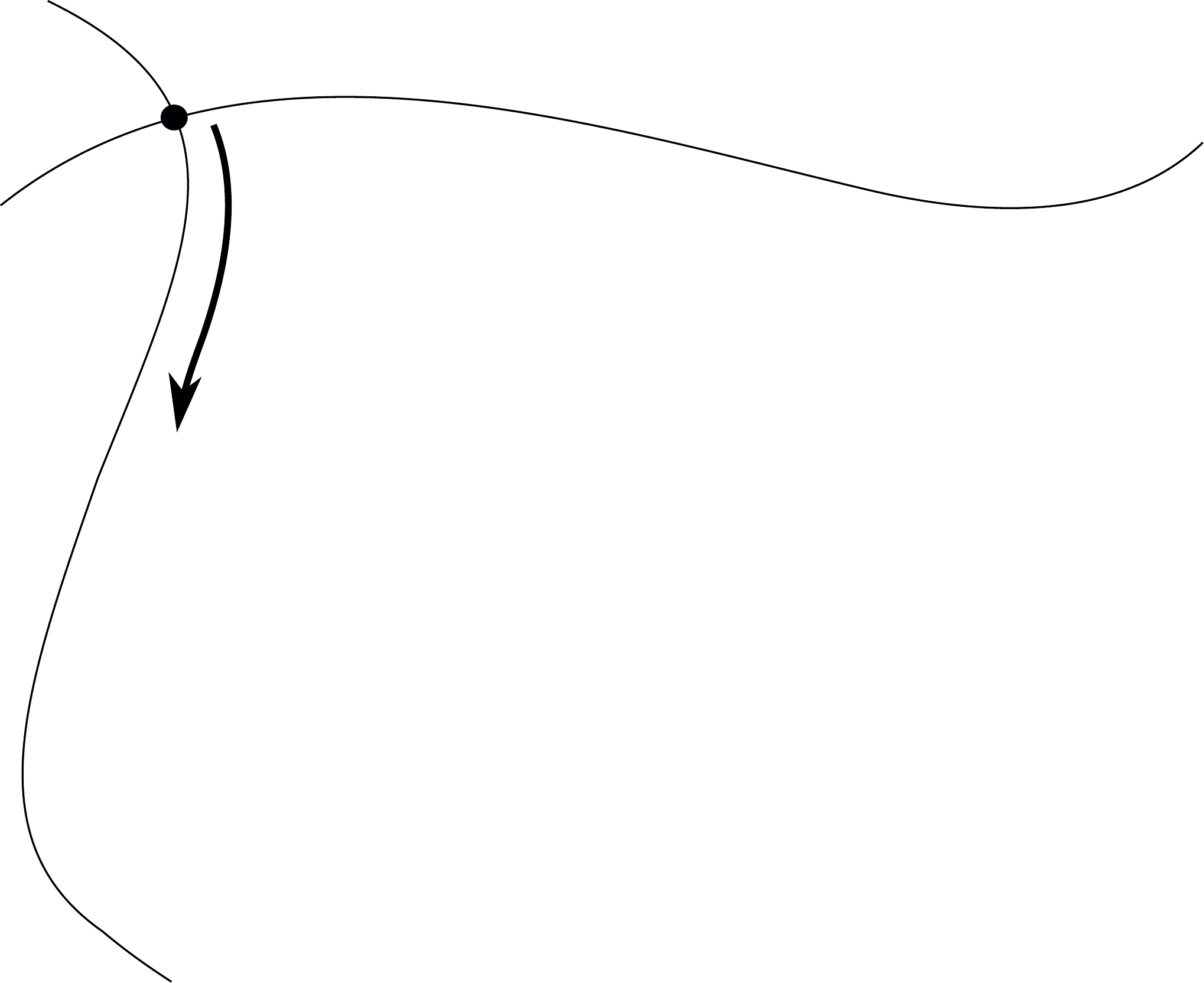}
\put (15, 74){$x$}
\put (65,72){$X$}

\put (45,63){$g(X)=i$}
\put (12,67){$y$}
\put (11,30){$g(Y)=g-i$}
\put (11,18){$Y$}
\end{overpic}
\end{center}
\end{figure}

We observe that $C\cdot D_\kappa=(2-2(g-i))c_{\min\{i,g-i\}}$. 
 
To calculate this intersection directly we must locate the limits of the type $\kappa$ in our test curve.  As $X$ and $Y$ are general curves and $x$ is a general point there are two possibilities. We have either $i$ points $p_j$ on $X$  or $i-1$ points $p_j$ on $X$.

In the first case we have $i$ points $p_j$ occur on $X$ indexed by the set $I$. 
In the $Y$-aspect we have
\begin{equation*}
(||I||-2i)y+\sum_{j\in I^C}k_jp_j\sim K_Y.
\end{equation*}
But as $|I^C|=g-i-2$ we only have solutions to the above if $||I||>2i$. In this case there are
\begin{equation*}
\dJ[g-i;||I||-2i,k_j\text{ for $j\in I^C$}]
\end{equation*}
solutions.
In the $X$-aspect we have
\begin{equation*}
(||I^C||-2(g-i))x+\sum_{j\in I}k_jp_j\sim K_X.
\end{equation*}
As $|I|=i$ and $||I||>2i$ we have $|I^C|=g-i-2$ and $||I^C||\leq 2g-2i-2$. By the Picard variety method (or De Jonquieres) there are 
\begin{equation*}
i!\prod_{i\in I}k_i^2
\end{equation*}
solutions. But we must discount for any solutions with $p_j=x$. For each $k_j>2(g-i)-||I^C||$ there are
\begin{equation*}
\dJ[i;k_i\in I-\{j\}]
\end{equation*}
 solutions each with order $k_j+||I^C||-2(g-i)+1$ by the order of the intersection of the components of the ramification loci of the relevant Jacobi map at these points.
 
 In the second case we have $i-1$ points $p_j$ occur on $X$ indexed by the set $I$. 
In the $Y$-aspect we have
\begin{equation*}
(||I||-2i)y+\sum_{j\in I^C}k_jp_j\sim K_Y.
\end{equation*}
 Now as $|I^C|=g-i-1$ we find by the Picard variety method when $||I||-2i\leq-2$ that there are 
 \begin{equation*}
 (g-i)!(||I||-2i)^2\prod_{j\in I^C}k_j^2
 \end{equation*}
 solutions. The case that $||I||-2i>-2$ will not be possible in the $X$-aspect. But we must discount for any solutions where $p_j=y$. This can occur only if $k_j>2i-||I||$ and in this case there are
 \begin{equation*}
 \dJ[g-i;k_j+(||I||-2i),k_i \text{ for $i\in I^C-\{j\}$}]
\end{equation*}
 such solutions each with order $k_j+(||I||-2i)+1$. In the $X$-aspect we have
 \begin{equation*}
(||I^C||-2(g-i))x+\sum_{j\in I}k_jp_j\sim K_X.
\end{equation*}
 for which there are 
\begin{equation*}
\dJ[i;k_j\text{ for $j\in I$}] 
\end{equation*}
 solutions when $||I^C||\geq2(g-i)$ and no solutions otherwise. There are no solutions with any $p_j=x$.

The last case that we must account for is if there is a zero sitting on a $\PP^1$-bridge between $x$ and $y$. First we provide a simple argument from the perspective of $g^r_d$'s in the case that $y$ is not a Weierstrass point. We insert a $\PP^1$-bridge glued at $s$ to $x$ and $t$ to $y$. If $y$ is not a Weierstrass point we know the vanishing sequence at $y$ to be $(i-1,i,...,\overline{2i-1},...,g+i-1)$ and at $x$ it is $(g-i-1,g-i,...,\overline{2(g-i)-1},...,2g-i-1)$ where we adopt the convention that the lined integer is omitted from the sequence. This gives us vanishing sequence at $t$ as $(g-i-1,g-i,...,\overline{2(g-i)-1},...,2g-i-1)$ and at $s$ as $(i-1,i,...,\overline{2i-1},...,g+i-1)$. In coordinates $[S;T]$ we let $s=[0;1]$ and $t=[1;0]$. We immediately observe that such a $g^{g-1}_{2g-2}$ must contain sections 
\begin{equation*}
S^{g+i-1}T^{g-i-1},S^{g+i-2}T^{g-i},..., \overline{S^{2i-1}T^{2(g-i)-1}},...,S^{i-1}T^{2g-i-1}.
\end{equation*}
But these sections are independent and hence form a basis of our $g^{g-1}_{2g-2}$. An isolated zero of order $k_m$ occurring at $[a;b]\ne s,t$ would be a section of the form
\begin{equation*}
S^nT^m(bS-aT)^{k_m}
\end{equation*}
with $n+m+k_m=2g-2$. But such a section would place $S^{2i-1}T^{2(g-i)-1}$ in our $g^{g-1}_{2g-2}$ unless $n>2i-1$ and $m>2(g-i)-1$ which is not possible, providing a contradiction unless $k_m=0$. 

The other case we need to consider is the case that $y$ is a normal Weierstrass point. We provide an argument based on the residues of the meromorphic differentials. See \cite{EisenbudHarrisWeierstrass} \S2 for discussion. Let $I$ index the $i$ points that sit across $X$ and the $\PP^1$-bridge and then $I^C$ indexes the points on $Y$. In this case as $Y$ is a general curve if we have 
\begin{equation*}
(||I||-2i)y+\sum_{j\in I^C}k_jp_j\sim K_Y.
\end{equation*}
on the $Y$-aspect and $y$ is also a Weierstrass point then these must be the same condition and we have $k_j=1$ for $j\in I^C$. This implies $||I^C||=g-i-2$ and $||I||=g+i$. If the isolated zero on the $\PP^1$-bridge has order $k_m$ then in the $X$-aspect we have
\begin{equation*}
(||I^C||-2(g-i)+k_m)x+\sum_{j\in I-\{m\}}k_jp_j\sim K_X.
\end{equation*}
As $x$ is a general point this only has solutions for $||I^C||-2(g-i)+k_m\geq 0$ which implies $k_m\geq g-i+2$. Finally we consider the meromorphic differential we are left with on the $\PP^1$-bridge. In the $\PP^1$-aspect we have
\begin{equation*}
-(g-i+2)t+k_mp+(g-i-k_m)s\sim K_{\PP^1}.
\end{equation*}
Now as the residues at $x$ and $y$ are zero we require also the residue at $s$ and $t$ to be zero. Placing $s$ at $0$, $t$ at $\infty$ and $p$ at $1$ we observe that locally our differential is
\begin{equation*}
\frac{(z-1)^{k_m}}{z^{k_m-(g-i)}}dz
\end{equation*}
which has zero residue only if $k_m\leq g-i$ thus providing a contradiction.

Putting the above together we have
\begin{eqnarray*}
&&C_i\cdot D_\kappa=\sum_{|I|=i,||I||\geq2i+1}\dJ[g-i;||I||-2i,k_j\text{ for $j\in I^C$}]\\
&&\left(i!\prod_{i\in I}k_i^2-\sum_{j\in I}(k_j+||I^C||-2(g-i)+1)\dJ[i;k_i\text{ for $i\in I-\{j\}$}]  \right)\\
&&+\sum_{|I|=i-1,||I||\leq 2i-2} \dJ[i;k_j\text{ for $j\in I$}]
\biggl((g-i)!(||I||-2i)^2\prod_{j\in I^C}k_j^2\\
&& -\sum_{\tiny{\begin{array}{cc}j\in I^C,\\k_j\geq 2i-||I||+1\end{array}}}(k_j+(||I||-2i)+1) \dJ[g-i;k_j+(||I||-2i),k_i \text{ for $i\in I^C-\{j\}$}]    \biggr).
\end{eqnarray*}

\subsection{The class of $D_\kappa$}\label{lab:formula}
The test curve results imply
\begin{eqnarray*}
c_1&=& \frac{-1}{2(g-2)}\left(  4(g-1)!\prod k_i^2+\sum_{k_i\geq 3}(k_i^2-k_i)\dJ[g-1;k_1,...,k_i-2,...,k_{g-2}]  \right),\\
c_{\min\{i,g-i\}}&=& \frac{-1}{2(g-i)-2}C_i\cdot D_\kappa,\hspace{0.5cm}\text{ where $C_i\cdot D_\kappa$ is given above},   \\
c_0&=&\frac{-1}{2(g-1)}\biggl( (g-1)!\prod k_i^2+\sum_{i=1}^{g-2}k_i\left(\sum_{j=1}^{k_i-1} \dJ[g-1;j-1,k_1,...,k_{i-1},k_{i+1},...,k_{g-2}]     \right)\\
&&+ \frac{1}{2(g-2)}\biggl(   4(g-1)!\prod k_i^2+\sum_{k_i\geq 3}(k_i^2-k_i)\dJ[g-1;k_1,...,k_i-2,...,k_{g-2}] \biggr)                                     \biggr),\\
c_\lambda&=& c_1-12c_0\\ 
&=&\frac{7-g}{2(g-1)(g-2)}\left(  4(g-1)!\prod k_i^2+\sum_{k_i\geq 3}(k_i^2-k_i)\dJ[g-1;k_1,...,k_i-2,...,k_{g-2}]  \right)\\
&&+\frac{6}{g-1}\left((g-1)!\prod k_i^2+\sum_{i=1}^{g-2}k_i\left(\sum_{j=1}^{k_i-1} \dJ[g-1;j-1,k_1,...,k_{i-1},k_{i+1},...,k_{g-2}]     \right)\right).
\end{eqnarray*}

\section{Components of $D_\kappa$}\label{sec:comp}

\subsection{$\kappa=(4)$ for $g=3$}
In this case Kontsevich and Zorich \cite{KontsevichZorich} showed that there are two connected components corresponding to the hyperelliptic component and the odd spin structure component. By our previous formula we have
\begin{equation*}
D_{(4)}=380\lambda-40\delta_0-100\delta_1.
\end{equation*}
Denoting the hyperelliptic component $D_{(4)}^\text{hyp}$ we observe that 
\begin{equation*}
C_1\cdot D_{(4)}^\text{hyp} =(3+5)\cdot 6=-2c_1
\end{equation*}
as $C_1$ intersects $D_{(4)}^\text{hyp}$ only when $y$ is a Weierstrass point in $Y$. When $y$ is a Weierstrass point in $Y$ then the possible limits of Weierstrass points on the nodal curve are the $3$ points on $X$ that are $2$-torsion to $x$ and the $5$ other Weierstrass points on $Y$. Similarly we observe
\begin{equation*}
B\cdot D_{(4)}^\text{hyp}=6+2=-4c_0+c_1
\end{equation*}
as $B$ intersects $D_{(4)}^\text{hyp}$ only when $x=y'$. In this case we have by the theory of admissible covers, the degenerate double cover has $6$ ramification points on the $C$ component (the genus $2$ curve) and two ramification points on a $\PP^1$ component that is attached to the $C$ component at $x$ and $y$ which are conjugate under this double cover. Hence as $A\cdot  D_{(4)}^\text{hyp}=0$ we have shown $D_{(4)}^\text{hyp}=8H$ where\begin{equation*}
H=9\lambda-\delta_0-3\delta_1
\end{equation*}
as expected as the general hyperelliptic curve of genus $g=3$ has $8$ branch points. Hence
\begin{equation*}
D_{(4)}^\text{odd}=308\lambda-32\delta_0-76\delta_1
\end{equation*}
which matches the calculation in \cite{Cukierman} \S5 of the closure of the loci of hyperflexes on plane quartics.

\subsection{$\kappa=(3,3)$ for $g=4$}
Kontsevich and Zorich \cite{KontsevichZorich} show that in this case there are two components corresponding to the hyperelliptic component and the non-hyperelliptic component. In our case the hyperelliptic component drops dimension when we pushdown from $\overline{\mathcal{M}}_{4,2}$ and hence we have calculated the class of the non-hyperelliptic component
\begin{equation*}
D_{(3,3)}=D_{(3,3)}^\text{non-hyp}.
\end{equation*}

\subsection{$\kappa=(4,2^{g-3})$ for $g\geq 4$}
In this case Kontsevich and Zorich \cite{KontsevichZorich} showed that there are two connected components corresponding to the odd and even spin structure components. Teixidor i Bigas \cite{Teixidor} calculated $D_\kappa^\text{even}$ and these divisors are the pushdown to $\Mgb$ of divisors calculated by Farkas and Verra \cite{Farkas}\cite{FarkasVerra} on Cornelia's compactified moduli space of spin curves. Our calculation will be from the point of view of degenerating abelian differentials and will employ the theory of admissible covers. The theory of limit linear series in the case that $r=1$ is known as the theory of admissible covers as degenerating $g_d^r$'s in this case can be considered as degenerating degree $d$ covers of $\PP^1$.

Both components have zero intersection with test curve $A$ and we are left to apply test curves $B$ and $C_i$ to each component.

\subsubsection{Test curve $B$} We have discussed the two different ways that spin structures manifest on a curve with a non-separating node. For the test curve $B$ we must decipher the spin structure of the different solutions.

Consider a smooth curve $Z$ such that
\begin{equation*}
4p_1+\sum_{i=2}^{g-2}2p_i\sim K_Z.
\end{equation*}
The associated theta characteristic is
\begin{equation*}
\eta\sim 2p_1+\sum_{i=2}^{g-2}p_i
\end{equation*}
and the spin structure parity is equal to the parity of $h^0(Z,\eta)$. Because $\eta$ is effective we know that $h^0(Z,\eta)\geq 1$. Further, if $h^0(Z,\eta)=2$ then the sections give a degree $g-1$ cover of $\PP^1$. The theory of admissible covers tells us how such a map degenerates when we degenerate to a nodal curve and hence in this case how the theta characteristic degenerates. There are two cases of admissible covers where the domain curve stably reduces to a smooth genus $g-1$ curve with points $x$ and $y$ identified that arise in our test curve $B$. First we have a degree $g-1$ cover of $\PP^1$ by $C$ where $x$ and $y$ lie above the same point in $\PP^1$ with a $\PP^1$ connecting these two points in a different component of the cover. All other conjugate points to $x$ and $y$ in the cover must have degree one $\PP^1$-tails attached as pictured. In general there may be ramification at $x$ and/or $y$, however in the test curve $B$ this case does not show up. 

\begin{figure}[htbp]
\begin{center}
\begin{overpic}[width=0.5\textwidth]{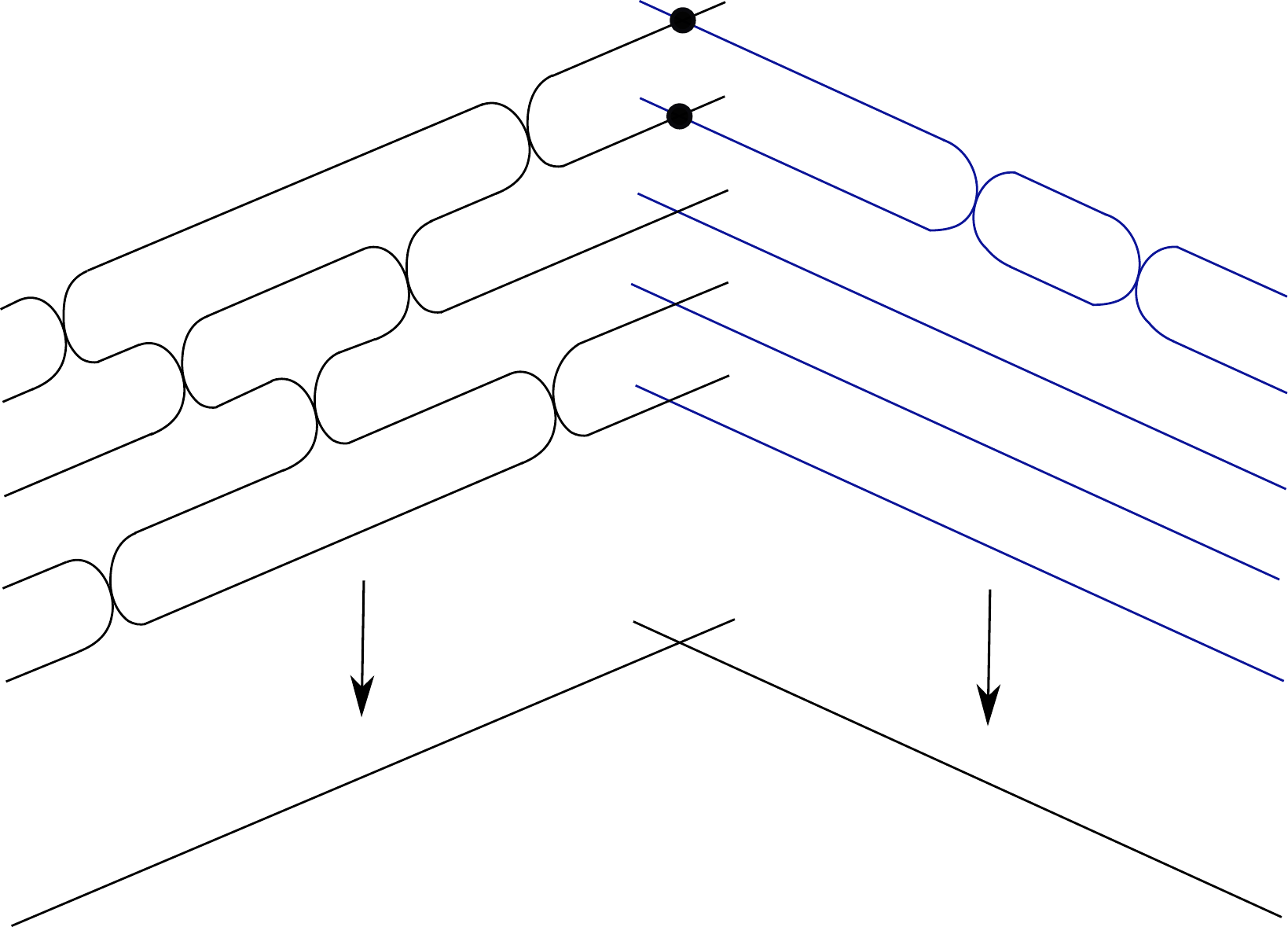}

\put (-8,31){$C$}
\put (-8,0){$\PP^1$}
\put (103,0){$\PP^1$}
\put (103,33){$\PP^1$}
\put (103,43){$\PP^1$}
\put (103,25){$\PP^1$}
\put (103,17){$\PP^1$}

\put (52,66){$x$}
\put (52,59){$y$}

\end{overpic}
\end{center}
\end{figure}

In fact, if we have a solution of the type 
\begin{equation*}
2x+2y+\sum_{i=2}^{g-2}2p_i\sim K_C+x+y,
\end{equation*}
then as $C$ is a general curve we have 
\begin{equation*}
h^0(C,{\eta})=h^0(C, x+y+\sum_{i=2}^{g-2}p_i )=h^0(C,K_C-\sum_{i=2}^{g-2}p_i)=2 
\end{equation*}
and this provides our cover of $\PP^1$ showing that this solution is indeed the limit of canonical divisors of type $\kappa$ with even spin structure on smooth curves.

The second type of nodal solution that occurs in our test curve is a solution of the type
\begin{equation*}
3x+y+\sum_{i=2}^{g-2}2p_i\sim K_C+x+y.
\end{equation*}
In this case we see that the spin structure has degenerated to a spin structure of the second kind on a nodal curve of this type with a square root $\eta$ of $K_C$ and $\OO(1)$ on a $\PP^1$-bridge between $x$ and $y$. If $h^0(C,\eta)=2$ then the theory of admissible covers dictates that the cover must be of the form pictured with $x$ and $y$ sitting above different points in the cover, but connected by $\PP^1$-bridges that would contract to form a node under stable reduction.
\vspace{0.5cm}
\begin{figure}[htbp]
\begin{center}
\begin{overpic}[width=0.8\textwidth]{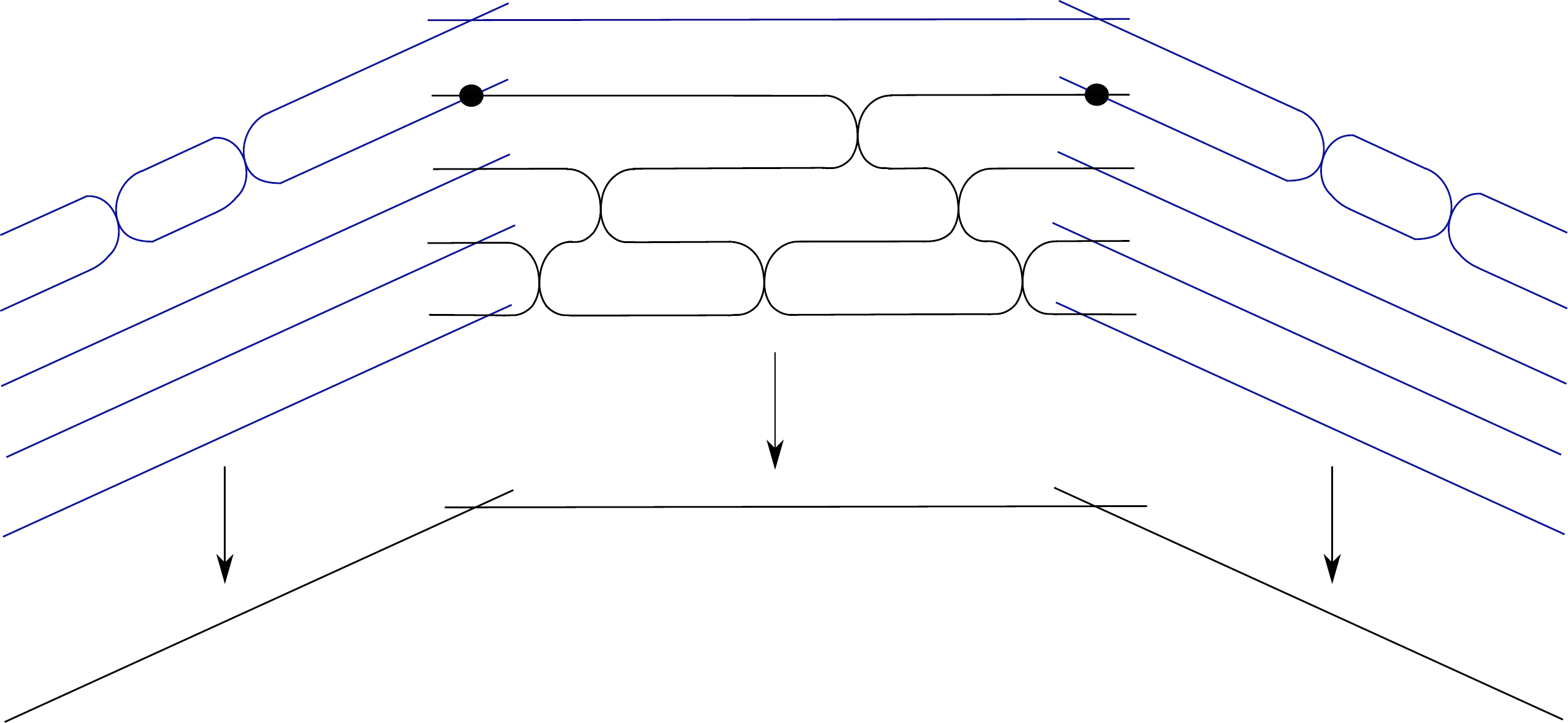}
\put (37,22){$C$}
\put (37,46){$\PP^1$}

\put (-5,0){$\PP^1$}
\put (-5,10){$\PP^1$}
\put (-5,15){$\PP^1$}
\put (-5,20){$\PP^1$}
\put (-5,27){$\PP^1$}

\put (102,0){$\PP^1$}
\put (102,10){$\PP^1$}
\put (102,15){$\PP^1$}
\put (102,20){$\PP^1$}
\put (102,27){$\PP^1$}

\put (69,37){$x$}
\put (29,37){$y$}

\end{overpic}
\end{center}
\end{figure}

However, in this case we see that $y$ forms part of a special section as
\begin{equation*}
y+\sum_{i=1}^{g-3}q_i\sim \eta
\end{equation*}
and hence
\begin{equation*}
2y+\sum_{i=1}^{g-3}2q_i\sim K_C.
\end{equation*}
But this is not possible in our case as we have chosen $y$ to be a general point. Hence the nodal solution of this type in our test curve are located in the odd spin structure component.

Hence solutions in our test curve with even spin structure arise as ramification points of a cover of the first type. Riemann-Hurwitz gives us that there are $4g-4$ simple ramification points for any such cover of $\PP^1$. This includes the $2$ ramification points on the $\PP^1$-bridge between $x$ and $y$ which become an order $2$ solution at the node. For a fixed $y$ in a general curve $C$ there are $\dJ[g-1;1,2^{g-3}]$ such $x$ and hence we obtain
\begin{equation*}
B\cdot D_\kappa^\text{even}=(4g-4)\dJ[g-1;1,2^{g-3}].
\end{equation*}
We are now presented with the question if any of these solutions are also limits of odd spin structure of the type we are interested in. These solutions are all sections of $H^0(C,\tilde{\eta})$ which descend to sections of $\eta$ on the nodal curve under the given gluing. Only the nodal section still descends when we change this gluing and hence it is only the nodal solution that is also in the odd spin structure component. Hence we obtain
\begin{equation*}
B\cdot D_\kappa^\text{odd}=(g-1)!4^2\cdot 2^{2(g-3)}+4(\dJ[g-1; 1;2^{g-3}]+\dJ[g-1;2^{g-2}])-(4g-4)\dJ[g-1;1,2^{g-3}].
\end{equation*} 
We now compute the relevant values of de Jonquieres' formula 
\begin{eqnarray*}
\dJ[g-1;2^{g-2}]
&&=(g-1)! \cdot 2^{g-2}\left(\sum_{j=0}^{g-3}\frac{(-1)^j}{1+j}\left(\begin{pmatrix}  g-2\\g-2-j\end{pmatrix}2^{g-2-j}   \right)   +\frac{(-1)^{g-2}}{g-1}\right)\\
&&=(g-1)!\cdot 2^{g-2}\left(\frac{(-1)^{g-1}}{g-1}+\frac{2^g-2}{2(g-1)}+\frac{(-1)^{g-2}}{g-1}   \right)\\
&&=(g-2)!\cdot 2^{g-2}(2^{g-1}-1)
\end{eqnarray*}
and
\begin{eqnarray*}
\dJ[g-1;1,2^{g-3}]
&&=(g-1)! \cdot 2^{g-3}\left(\sum_{j=0}^{g-3}\frac{(-1)^j}{1+j}\left(\begin{pmatrix}  g-3\\g-2-j\end{pmatrix}2^{g-2-j} +\begin{pmatrix}  g-3\\g-3-j\end{pmatrix}2^{g-3-j}   \right)   +\frac{(-1)^{g-2}}{g-1}\right)\\
&&=(g-1)!\cdot 2^{g-3}\left(\frac{4(g-2)(-1)^{g-2}+4+(8-4g)(-1)^g+(g-3)2^g }{4(g-2)(g-1)}  \right)\\
&&=(g-3)!\cdot 2^{g-3}((g-3)2^{g-2}  +1 ).
\end{eqnarray*}
This gives
\begin{eqnarray*}
B\cdot D_\kappa^\text{even}&=&  (4g-4)(g-3)!\cdot 2^{g-3}((g-3)2^{g-2}  +1 )      \\
B\cdot D_\kappa^\text{odd}&=& (g-2)!\cdot 2^{g-3}((g+5)2^g-12).       \\
\end{eqnarray*}

\subsubsection{$C_1$} The spin structure on a nodal curve of compact type  comes from the theta characteristics on each component $(\eta_Y,\eta_X)$. There are two types limits that we must consider here. We either have a zero of order $4$ on the $X$ component or no zero on the $X$ component. In the first situation we have
\begin{eqnarray*}
&&\eta_Y\sim y+\sum_{i=2}^{g-2}p_i,\\
&&\eta_X\sim 2p_1-2x.
\end{eqnarray*}
We observe $h^0(Y,\eta_Y)=1$ as $Y$ is a general curve. We also observe that $h^0(X,\eta_X)=1$ for the $3$ solutions for $p_1$ which are $2$-torsion points to $x$ and zero for the other $12$ solutions. Hence we have a contribution of $3\dJ[g-1;2^{g-2}]$ to the intersection with the even component and $12\dJ[g-1;2^{g-2}]$ to the intersection with the odd component.

For the solutions with no zero on  the $X$ component we have
\begin{eqnarray*}
&&\eta_Y\sim -y+2p_1+\sum_{i=2}^{g-2}p_i,\\
&&\eta_X\sim \OO_X.
\end{eqnarray*}
We know $h^0(X,\eta_X)=1$. As $\eta_Y^{\otimes 2}=K_Y$ we have
\begin{equation*}
h^0(Y,\eta_Y)=h^0(Y,K_Y-\eta_Y)=h^0\biggl(Y,K_Y+y-2p_1-\sum_{i=2}^{g-2}p_i\biggr). 
\end{equation*}
But as we know that $y$ is a base point of $K_Y+y$ we have $h^0(Y,K_Y+y-2p_1-\sum_{i=2}^{g-2}p_i)=h^0(Y,K_Y-2p_1-\sum_{i=2}^{g-2}p_i)$.

To find the solutions where $h^0(Y,\eta_Y)=1$ we will approach this problem by a reverse construction. If $h^0(Y,\eta_Y)=1$ we have
\begin{equation*}
2p_1+\sum_{i=2}^{g-2}p_i+\sum_{i=1}^{g-3}q_i\sim K_C  
\end{equation*}
for some $q_i$. But we already know
\begin{equation*}
4p_1+\sum_{i=2}^{g-2}2p_i -2y\sim K_C.
\end{equation*}
Hence we obtain that $h^0(Y,\eta_Y)= 1$ if and only if there exist $q_i$ such that 
\begin{equation*}
2y+\sum_{i=1}^{g-3}2q_i\sim K_C.
\end{equation*}
The $p_i$ are then the sets of solutions such that
\begin{equation*}
2p_1+\sum_{i=2}^{g-2}p_i+\sum_{i=1}^{g-3}q_i\sim K_C
\end{equation*}
for each set of $q_i$.

We first observe that there are $\dJ[g-1;2^{g-2}]/(g-3)!$ solutions for
\begin{equation*}
2y+\sum_{i=1}^{g-3}2q_i\sim K_Y
\end{equation*} 
where we are considering the $q_i$ to be unordered. For each of these solutions we observe that there are $(4g-6)(g-3)!$ solutions to
\begin{equation*}
2p_1+\sum_{i=2}^{g-2}p_i+\sum_{i=1}^{g-3}q_i\sim K_Y
\end{equation*}
which are just the ramification points of $|K_Y-\sum_{i=1}^{g-3}q_i|$. Now we just need to omit the solutions where $y=p_1$ of which there are $\dJ[g-1;2^{g-2}]$. This gives a contribution of $(4g-7)\dJ[g-1;2^{g-2}]$ to the intersection with the even component. This leaves a contribution of
\begin{equation*}
4(g-1)!4^2\cdot 2^{2(g-3)}-3\dJ[g-1;2^{g-2}]-(4g-7)\dJ[g-1;2^{g-2}]=2^{g}(g-1)!(2^g-2^{g-1}+1)
\end{equation*}
to the intersection with the odd component. Putting this together we have
\begin{eqnarray*}
C_1\cdot D_\kappa^\text{even}&=&  (4g-7)\dJ[g-1;2^{g-2}]      +3\dJ[g-1;2^{g-2}]     \\
&=& (g-1)!\cdot 2^{g}(2^{g-1}-1) ,                \\
C_1\cdot D_\kappa^\text{odd}&=&  2^{g}(g-1)!(2^g-2^{g-1}+1)+12\dJ[g-1;2^{g-2}]         \\
&=&2^{g-1}(g-2)!((g+2)2^g+2g-8).
\end{eqnarray*} 
This gives
\begin{eqnarray*}
c_1^\text{even}&=&-(g-1)(g-3)!\cdot 2^{g-1}(2^{g-1}-1)   ,             \\
c_1^\text{odd}&=&-2^{g-2}(g-3)!((g+2)2^g+2g-8),
\end{eqnarray*}
and hence
\begin{eqnarray*}
c_0^\text{even}&=&   -2^{2g-4}(g-1)(g-3)!                 \\
c_0^\text{odd}&=& -2^{g-4}((g+6)2^g-8)(g-3)!
\end{eqnarray*}
and
\begin{eqnarray*}
c_\lambda^\text{even}&=&   2^{g-1}(2^g+1)(g-1)(g-3)!                 \\
c_\lambda^\text{odd}&=& 2^{g-1}(2^g-1)(g+8)(g-3)!.
\end{eqnarray*}

\subsubsection{Test curve $C_i$}
The spin structure on a nodal curve of compact type  comes from the theta characteristics on each component $(\eta_Y,\eta_X)$. There are two types limits that we must consider here. We either have $i$ zeros of the form $\{4,2^{i-1}\}$ on the $X$ component or we have $i-1$ zeros of the form $\{2^{i-1}\}$ on the $X$ component. In the first situation we have
\begin{eqnarray*}
&&\eta_Y\sim y+\sum_{j=i+1}^{g-2}p_j,\\
&&\eta_X\sim -2x+2p_1+\sum_{j=2}^{i}p_j.
\end{eqnarray*}
We have $h^0(Y,\eta_Y)=1$ as $Y$ is a general curve. As $X$ is a general curve and $x$ is a fixed general point, we have  $h^0(\eta_X)$ is ether $1$ or $0$. If $h^0(\eta_X)=1$ we have
\begin{equation*}
\sum_{j=1}^{i-1}q_j\sim \eta_X
\end{equation*}
and hence we have 
\begin{equation*}
\sum_{j=1}^{i-1}2q_j\sim K_X.
\end{equation*}
In fact, for any set of $q_i$ satisfying this equation we have $p_i$ arising as the ramification of the series $|\sum_{j=1}^{i-1}q_j+2x|=|\eta_X+2x|$. There are $\dJ[i;2^{i-1}]/(i-1)!$ such sets of $q_j$ (where we are not allowing for the order of the $q_j$) and each set of $q_j$ has $4i$ ramification points. We just need to discount for the solution where $p_1=x$, which is unique for each choice of $q_j$ up to the ordering of the $p_j$ for $j=2,...,i$. This gives a contribution of 
\begin{equation*}
\frac{(4i-1)\dJ[i;2^{i-1}]\dJ[g-i;2^{g-i-1}](g-i-1)(g-3)!}{(i-1)!(g-i-1)!}
\end{equation*}
to the intersection with the even component.

In the second case if we have zeros of the form $\{2^{i-1}\}$ on the $X$ component we have
\begin{eqnarray*}
&&\eta_Y\sim -y+2p_1+\sum_{j=i+1}^{g-2}p_j,\\
&&\eta_X\sim \sum_{j=2}^{i}p_j.
\end{eqnarray*}
We observe $h^0(\eta_X)=1$ as $X$ is a general curve. If $h^0(\eta_Y)=1$ then we have
\begin{equation*}
\sum_{j=1}^{g-i-1}q_j\sim\eta_Y
\end{equation*}
and hence
\begin{equation*}
\sum_{j=1}^{g-i-1}2q_j\sim K_Y.
\end{equation*}
Again, for any such set of $q_j$ we can obtain the $p_j$ and $y$ and we have by the Picard variety method that there is a contribution of
\begin{equation*}
\frac{\dJ[i;2^{i-1}]}{(i-1)!}\left(   \frac{\dJ[g-i;2^{g-i-1}]}{(g-i-1)!} \left((-1)^2 2^2(g-i)(g-i-1) -3(g-i-1) \right)  \right)(g-3)!
\end{equation*}
to the intersection with the even component. The correction term comes from enumerating the solutions with $y=p_1$ which have an order of $3$ coming from our knowledge of their relation in the canonical embedding of $X$.

Putting this together we have
\begin{eqnarray*}
C_i\cdot D_\kappa^\text{even}&=&   \frac{\dJ[i;2^{i-1}]}{(i-1)!}  \frac{\dJ[g-i;2^{g-i-1}]}{(g-i-1)!} (g-3)!((4i-1)+4(g-i)-3)(g-i-1)    \\
&=&2^{g}(2^{g-i}-1)(2^i-1)(g-1)(g-i-1)(g-3)! \\
C_i\cdot D_\kappa^\text{odd}&=&   2^{g-i}((2^i-1)(2^g+2^i)g+2(2^g-2^{2i})i)(g-i-1)(g-3)! 
\end{eqnarray*}
which gives

\begin{eqnarray*}
c_i^\text{even}&=&   -2^{g-1}(2^{g-i}-1)(2^i-1)(g-1)(g-3)!                 \\
c_i^\text{odd}&=&       -2^{g-i-1}((2^i-1)(2^g+2^i)g+2(2^g-2^{2i})i)(g-3)!      . 
\end{eqnarray*}

\section{The slope of divisors $D_\kappa$}
The slope of an effective divisor $D$ on $\Mgb$ is defined as 
\begin{equation*}
s(D)=\frac{c_\lambda}{\min\{-c_i\}}
\end{equation*}
and it is well known that $s(D)<\infty$ for any $D$ which is the closure of an effective divisor on $\Mg$ (\cite{HarrisMorrisonSlope}). In all known examples, the slope is equal to 
\begin{equation*}
s_0(D)=\frac{c_\lambda}{-c_0}
\end{equation*} 
and Farkas and Popa \cite{FarkasPopa} showed this holds for $g\leq 23$ and conjectured that this is always true. We observe that by the formula for the class of $D_\kappa$ for $\kappa=(k_1,...,k_{g-2})$ given in \S\ref{lab:formula} we have
\begin{equation*}
s_0(D_{\kappa})=12-\frac{2g-2}{1+2(g-2)(B/A)}
\end{equation*}
where
\begin{eqnarray*}
A&=&4(g-1)!\prod k_i^2+\sum_{k_i\geq 3}(k_i^2-k_i)\dJ[g-1;k_1,...,k_i-2,...,k_{g-2}],\\
B&=&(g-1)!\prod k_i^2+\sum_{i=1}^{g-2}k_i\left(\sum_{j=1}^{k_i-1} \dJ[g-1;j-1,k_1,...,k_{i-1},k_{i+1},...,k_{g-2}]  \right),
\end{eqnarray*}
and de Jonquieres' formula is provided in \S\ref{sec:dJ}. If $\kappa=(4,2^{g-3})$ then $D_\kappa$ has two components and we have
\begin{eqnarray*}
s_0(D_{(4,2^{g-3})}^\text{even})&=&  8+\frac{1}{2^{g-3}} ,         \\
s_0(D_{(4,2^{g-3})}^\text{odd})&=&       8+\frac{ (2^{g+1}-g)  }{ 2^{g-3}(g+6)   -1}     .
\end{eqnarray*}
Computations in low genus suggest that for fixed genus $s_0(D_\kappa)$ for general $\kappa$ is bounded below by $s(D_{(4,2^{g-3})}^\text{even})=s_0(D_{(4,2^{g-3})}^\text{even})$ which is asymptotically $8$ and bounded above by $s(D_{(g+1,1^{g-3})})=s_0(D_{(g+1,1^{g-3})})$ originally calculated by Diaz \cite{Diaz}
\begin{equation*}
s(D_{(g+1,1^{g-3})})=\frac{3(3g^2+3g+2)}{g(g+1)}=9+\frac{2}{g(g+1)},
\end{equation*}
which is asymptotically $9$.

It would be interesting to understand the geometric significance of this range of the asymptotic $s_0$ and conjectured slope. Also of interest is the birational significance of these divisors $D_\kappa$, including their positions in the pseudo-effective cone $\overline{\text{Eff}}(\Mgb)$, and the geometric significance of the birational model associated to each $D_\kappa$.

\bibliographystyle{plain}
\bibliography{base}
\end{document}